\documentstyle[amssymb,12pt]{article}

\def\fsimeq{\displaystyle{\buildrel \simeq\over\longrightarrow}}
\def\fpr{\displaystyle{\buildrel {p_R}\over \longrightarrow}}
\def\fpl{\displaystyle{\buildrel {p_L}\over \longrightarrow}}
\def\fc{\displaystyle{\buildrel \subset\over \longrightarrow}}

\def\Talpha#1{\vbox{\ialign{##\crcr
    $\alpha$\crcr\noalign{\kern2pt\nointerlineskip}
           $\hfil\displaystyle{#1}\hfil$\crcr}}} \def\Out{\mbox{Out}}

\def\Int{\mbox{Int}}

\def\Onabla#1{\vbox{\ialign{##\crcr
$\,\scriptstyle{0}$\crcr\noalign{\kern2pt\nointerlineskip}
           $\hfil\displaystyle{#1}\hfil$\crcr}}}
\def\ker{{\rm ker}}

\def\cala{{\cal A}}
\def\calb{{\cal B}}
\def\call{{\cal L}}
\def\calm{{\cal M}}
\def\caln{{\cal N}}

\def\otimesinf{\mathop{\otimes}}

\def\bbbone{\mbox{\rm 1\hspace
{-.6em} l}} \def\c{\Bbb C} \def\hom{{\mbox{Hom}}}
 \def\der{{\mbox{\scriptsize Der}}}
\def\gder{{\mbox{Der}}}

 \newtheorem{theorem}{THEOREM}
\newtheorem{lemma}{LEMMA}

\begin{document}

\baselineskip=0.7cm
\begin{center} {\Large\bf ON THE FIRST ORDER OPERATORS\\
IN BIMODULES}
\end{center}
\vspace{0.75cm}

\begin{center}
Michel DUBOIS-VIOLETTE \\
and\\
Thierry MASSON\\

\vspace{0.3cm}
{\small
Laboratoire de Physique Th\'eorique et Hautes
Energies\footnote{Laboratoire associ\'e au Centre National de la
Recherche Scientifique - URA D0063}\\ Universit\'e Paris XI, B\^atiment
211\\ 91 405 Orsay Cedex, France\\ flad@qcd.th.u-psud.fr \&
masson@qcd.th.u-psud.fr}

\end{center} \vspace{1cm}

\begin{center} July 17, 1995 \end{center}

\vspace {1cm}

\noindent L.P.T.H.E.-ORSAY 95/56

\newpage \begin{abstract}

We analyse the structure of the first order operators in bimodules introduced
by A. Connes. We apply this analysis to the theory of connections on bimodules
generalizing thereby several proposals.

\end{abstract}

\section{Introduction}

It is well known that the notion of linear partial differential operator admits
the following algebraic formulation. Let $\cala$ be a commutative associative
algebra over $\Bbb C$ and let $\calm$ and $\caln$ be two $\cala$-modules. A
linear mapping $D$ of $\calm$  into $\caln$ is called an operator of order 0 of
$\calm$ into $\caln$ if it is a $\cala$-module homomorphism, i.e. if $[D,f]=0$,
$\forall f\in \cala$, ($f\in \cala$ being identified with the multiplication by
$f$ in $\calm$ and in $\caln$). Then, one defines inductively the operators of
order $k\in \Bbb N$: A linear mapping $D$ of $\calm$ into $\caln$ is an
operator of order $k+1$ of $\calm$ into $\caln$ if $[D,f]$ is an operator of
order $k$ of $\calm$ into $\caln$, $\forall f\in \cala$. For bad choices of
$\cala$, it may happen that this notion is not very appealing, however it makes
sense and when $\cala$ is the algebra of smooth functions on a smooth manifold
$X$ and when $\calm$ and $\caln$ are the modules of smooth sections of two
smooth complex vector bundles $E$ and $F$ over $X$ then this notion is just the
usual notion of linear partial differential operators of order $k\in \Bbb N$ of
$E$ into $F$. The troubles start if one tries to replace $\cala$ by a
noncommutative associative algebra and $\calm$ and $\caln$ by left (or right)
$\cala$-module. Then the definition breaks down because multiplications by
elements $f$ and $g$ of $\cala$ do not commute. However, it was noticed by A.
Connes [3] that if one works with bimodules, there is a natural noncommutative
generalization of first order operators since the multiplications on the right
and on the left commute. Moreover, in [3] it was pointed out that the first
order condition for the generalized Dirac operator plays an important role in
the connection between the noncommutative Poincar\'e duality and cyclic
cohomology. On the other hand, it was observed at several places, e.g. [6][7],
that the natural generalization of the notion of module over a commutative
algebra is not necessarily the notion of left (or right) $\cala$-module over a
noncommutative algebra $\cala$ but can be a notion of bimodule over $\cala$.
Furthermore it was claimed in [7] that the latter point of view is unavoidable
if one wants to discuss reality conditions and if one takes  the Jordan
algebras of hermitian elements of complex $\ast$-algebras as the noncommutative
analogue of algebras of real functions (e.g. as in quantum theory).\\
In section 2 we recall the definition of first order operators in a form which
is convenient for our purposes.\\
In section 3 we give some basic general examples.\\
In section 4 we establish the general structure of first order operators and
describe the appropriate notion of symbol.\\
In section 5 we apply our result to the theory of connections on bimodules and
show the relation between the symbols and previously introduced twistings
(generalized transpositions) in the case of noncommutative generalizations of
linear connections [9].

\section{First order operators in bimodules}

In this paper, $\cala$ and $\calb$ are unital associative algebras over $\Bbb
C$, their units $\bbbone_{\cala}$ and $\bbbone_{\calb}$ will be simply denoted
by $\bbbone$ when no confusion arises. Let $\calm$ and $\caln$ be two
$(\cala,\calb)$-bimodules, (i.e. $\cala\otimes\calb^{op}$-modules), and let
$\call(\calm,\caln)$ be the vector space of all linear mappings of $\calm$ into
$\caln$. Among the elements $\call(\calm,\caln)$, one can distinguish several
natural subclasses. The most natural subspace of $\call(\calm,\caln)$ is the
space $\hom^{\calb}_{\cala}(\calm,\caln)$ of all $(\cala,\calb)$-bimodule
homomorphisms of $\calm$ into $\caln$. Other natural subspaces are the space
$\hom^{\calb}(\calm,\caln)$ of all right $\calb$-module homomorphisms of
$\calm$ into $\caln$ and the space $\hom_{\cala}(\calm,\caln)$ of all left
$\cala$-module homomorphisms of $\calm$ into $\caln$. However, from the point
of view of the $(\cala,\calb)$-bimodule structure, it was pointed out in [3]
that there is a more symmetrical subspace of $\call(\calm,\caln)$ which
contains both $\hom^{\calb}(\calm,\caln)$ and $\hom_{\cala}(\calm,\caln)$ which
we now describe. One has the following lemma.

\begin{lemma}
The following conditions a) and b) are equivalent for an element $D$ of
$\call(\calm,\caln)$.\\
a) for any $f\in\cala$, $m\mapsto D(fm)-fD(m)$ is a right $\calb$-module
homomorphism.\\
b) for any $g\in \calb$, $m\mapsto D(mg)-D(m)g$ is a left $\cala$-module
homomorphism.
\end{lemma}

\noindent {\bf Proof} [3]. Let $L_f$ be the left multiplication by $f\in \cala$
and let $R_g$ be the right multiplication by $g\in \calb$ in $\calm$ and
$\caln$. The condition a) reads $[[D,L_f],R_g]=0,\  \forall f\in \cala, \forall
g\in \calb$, and the condition b) reads $[[D,R_g],L_f]=0,\ \forall f\in \cala,
\forall g\in \calb$. On the other hand one has $[L_f,R_g]=0$ which implies that
$[[D,L_f],R_g]=[[D,R_g],L_f]$. $\square$\\

\noindent An element $D$ of $\call(\calm,\caln)$ which satisfies the above
equivalent conditions a) and b) is called [3] {\sl a first order operator or an
operator of order} 1 {\sl of} $\calm$ {\sl into} $\caln$. The set of all first
order operators of $\calm$ into $\caln$ is a subspace of $\call(\calm,\caln)$
which will be denoted by $\call_1(\calm,\caln)$:\\
$\call_1(\calm,\caln)=\{D\in \call(\calm,\caln) \vert [[D,L_f],R_g]=0$,
$\forall f\in \cala, \forall  g\in \calb\}$.\\
This terminology is of course suggested by the fact that when $\cala$ and
$\calb$ coincide with the algebra $C^\infty(X)$ of smooth functions on a smooth
manifold $X$ and when $\calm$ and $\caln$ are the smooth sections of smooth
vector bundles $E$ and $F$ over $X$, then $\call_1(\calm,\caln)$ is just the
space of ordinary first order differential operators of $E$ into $F$.\\

\noindent{\bf Remark}. One has $D\in \hom_{\cala}(\calm,\caln)\Leftrightarrow
[D,L_f]=0$ $\forall f\in \cala$ and $D\in
\hom^{\calb}(\calm,\caln)\Leftrightarrow [D,R_g]=0$ $\forall g\in \calb$.
Therefore, $\hom_{\cala}(\calm,\caln)$ and\linebreak[4]
$\hom^{\calb}(\calm,\caln)$ are subspaces of $\call_1(\calm,\caln)$. When $f$
runs over $\cala$, $[D,L_f]$ are the obstructions for $D$ to be a left
$\cala$-module homomorphism and thus, the condition a) of lemma 1 means that
theses obstructions are right $\calb$-module homomorphisms. The condition b)
can be formulated similarily by exchange of left and right.

\section{Examples}

\subsection{The structure of $\call_1(\cala,\caln)$ for a
$(\cala,\cala)$-bimodule $\caln$}

Recall that a {\sl derivation of} $\cala$ {\sl into} a bimodule $\caln$ over
$\cala$ is a linear mapping $\delta:\cala\rightarrow \caln$ satisfying
$\delta(fg)=\delta(f)g+f\delta(g),\ \forall f,g\in \cala$. The space of all
derivations of $\cala$ into $\caln$ is denoted by $\gder(\cala,\caln)$. For
each element $n\in\caln$ one defines a derivation $ad(n)\in \gder(\cala,\caln)$
by $ad(n)(f)=nf-fn,\ \forall f\in\cala$. The subspace $ad(\caln)$ of
$\gder(\cala,\caln)$ is denoted by $\Int(\cala,\caln)$ and its elements are
called {\sl inner (or interior) derivations of} $\cala$ {\sl into} $\caln$. By
the very definition, a $D\in \call(\cala,\caln)$ is a first order operator,
i.e. $D\in \call_1(\cala,\caln)$, if and only if one has:
\[
D(fg)=D(f)g + fD(g) - fD(\bbbone)g,\ \forall f,g\in \cala.
\]
It follows that the derivations of $\cala$ into $\caln$ are exactly the first
order operators of $\cala$ into $\caln$ which vanish on the unit $\bbbone$ of
$\cala$, i.e. one has:
\[
\gder(\cala,\caln)=\{D\in \call_1(\cala,\caln)\vert D(\bbbone)=0\}.
\]
One defines two projections $p_R$ and $p_L$ of $\call_1(\cala,\caln)$ onto
$\gder(\cala,\caln)$ by setting $p_R(D)(f)=D(f)-D(\bbbone)f$ and
$p_L(D)(f)=D(f)-fD(\bbbone)$. Notice that $p_L(D)-p_R(D)=ad(D(\bbbone))$ so
$Im(p_L-p_R)=\Int(\cala,\caln)$. It is clear that $D\mapsto
(p_L(D),D(\bbbone))$ and $D\mapsto (p_R(D),D(\bbbone))$ are both isomorphisms
of $\call_1(\cala,\caln)$ onto $\gder(\cala,\caln)\oplus\caln$. Thus one has
\[
\call_1(\cala,\caln)\simeq \gder(\cala,\caln)\oplus \caln.
\]
In fact,  one has $\ker(p_R)=\hom^{\cala}(\cala,\caln)$,
$\ker(p_L)=\hom_{\cala}(\cala,\caln)$ and the mapping $D\mapsto D(\bbbone)$ of
$\call(\cala,\caln)$ into $\caln$ induces isomorphisms\\
$\alpha_R:\hom^{\cala}(\cala,\caln) \fsimeq \caln$ and
$\alpha_L:\hom_\cala(\cala,\caln)\fsimeq\caln$, so the two above isomorphisms
of $\call_1(\cala,\caln)$ onto $\gder(\cala,\caln)\oplus \caln$ correspond to
the canonical splitting of the exact sequences of vector spaces
\[
0\rightarrow \hom^{\cala}(\cala,\caln)\fc \call_1(\cala,\caln)\fpr
\gder(\cala,\caln)\rightarrow 0
\]
and
\[
0 \rightarrow \hom_{\cala}(\cala,\caln)\fc \call_1(\cala,\caln)\fpl
\gder(\cala,\caln)\rightarrow 0
\]
associated to the inclusion $\gder(\cala,\caln)\subset \call_1(\cala,\caln)$.

Notice finally that one has in $\call_1(\cala,\caln)$
\[
\Int(\cala,\caln)=\gder(\cala,\caln)\cap
(\hom^{\cala}(\cala,\caln)+\hom_{\cala}(\cala,\caln)).
\]

\subsection{First order operators of $\cala$ into itself}

In the case where $\caln=\cala$, the space $\call_1(\cala,\cala),
\gder(\cala,\cala)$ and $\Int(\cala,\cala)$ will be simply denoted by
$\call_1(\cala),\gder(\cala)$ and $\Int(\cala)$. If $D_1$ and $D_2$ are
elements of $\call_1(\cala)$, it is easy to see that $[D_1,D_2]=D_1\circ
D_2-D_2\circ D_1$ is again an element of $\call_1(\cala)$. Thus
$\call_1(\cala)$ is a Lie algebra; it is also, in an obvious way, a module over
the center $Z(\cala)$ of $\cala$. The subspace $\gder(\cala)$ is a Lie
subalgebra and a $Z(\cala)$-submodule of $\call_1(\cala)$ and $\Int(\cala)$ is
a Lie ideal of $\gder(\cala)$ and also a $Z(\cala)$-submodule. The quotient
$\Out(\cala)=\gder(\cala)/\Int(\cala)$ is a Lie algebra and a
$Z(\cala)$-module. The spaces $\hom^{\cala}(\cala,\cala)$ and
$\hom_{\cala}(\cala,\cala)$ are both Lie ideals  and $Z(\cala)$-submodules of
$\call_1(\cala)$ and the exact sequences
\[
0\rightarrow \hom^{\cala}(\cala,\cala)\rightarrow \call_1(\cala) \fpr
\gder(\cala)\rightarrow 0
\]
and
\[
0\rightarrow \hom_{\cala}(\cala,\cala)\rightarrow \call_1(\cala) \fpl
\gder(\cala)\rightarrow 0
\]
are now exact sequences of Lie algebras and of $Z(\cala)$-modules, whereas the
corresponding isomorphisms $\call_1(\cala)\simeq \gder(\cala)\oplus \cala$ are
isomorphisms of Lie algebras and of $Z(\cala)$-modules. One has the isomorphism
of Lie algebra and of $Z(\cala)$-module
\[
\Out(\cala)\simeq
\call_1(\cala)/(\hom^{\cala}(\cala,\cala)+\hom_{\cala}(\cala,\cala)).
\]

\subsection{First order operators associated with derivations}

Let $\Omega^1_L$ be a bimodule over $\cala$, let $\Omega^1_R$ be a bimodule
over $\calb$ and let $d_L$ be a derivation of $\cala$ into $\Omega^1_L$ and
$d_R$ be a derivation of $\calb$ into $\Omega^1_R$. Let $\calm$ and $\caln$ be
two $(\cala,\calb)$-bimodules and let $D$ be a linear mapping of $\calm$ into
$\caln$.

\begin{enumerate}
\item
Assume that there is a $(\cala,\calb)$-bimodule homomorphism $\sigma_L$ of
$\Omega^1_L\otimes_{\cala}\calm$ into $\caln$ such that
$D(fm)=fD(m)+\sigma_L(d_L(f)\otimes m),\ \forall m\in \calm$ and $\forall f\in
\cala$, then $D$ is a first order operator.
\item
Assume that there is a $(\cala,\calb)$-bimodule homomorphism $\sigma_R$ of
$\calm\otimesinf_{\calb}\Omega^1_R$ into $\caln$ such that
$D(mg)=D(m)g+\sigma_R(m\otimes d_R(g))$, $\forall m\in \calm$ and $\forall g\in
\calb$, then $D$ is a first order operator.
\end{enumerate}

Thus any of the two above conditions implies that $D$ is of first order. We
shall now show that, conversely, if $D$ is a first order operator of $\calm$
into $\caln$, these conditions are satisfied with an appropriate choice of the
$(d_L,\Omega_L)$ and $(d_R,\Omega_R)$.

\section{General structure of first order operators}

Recall that in the category of derivations of $\cala$ into the bimodules over
$\cala$ there is an initial object, $d_u:\cala\rightarrow \Omega^1_u(\cala)$,
which is obtained by the following standard construction [2], [1]. The bimodule
$\Omega^1_u(\cala)$ is the kernel of the multiplication $m:\cala\otimes
\cala\rightarrow \cala$, $(m(f\otimes g)=fg)$, and the derivation
$d_u:\cala\rightarrow \Omega^1_u(\cala)$ is defined by $d_u(f)=\bbbone \otimes
f - f\otimes \bbbone$. This derivation has the following universal property :
{\sl For any derivation} $\delta:\cala\rightarrow \calm$ {\sl of} $\cala$ {\sl
into a bimodule} $\calm$ {\sl over} $\cala$, {\sl there is a unique bimodule
homomorphism} $i_\delta$ {\sl of} $\Omega^1_u(\cala)$ {\sl into} $\calm$ {\sl
such that} $\delta=i_\delta \circ d_u$. As left (resp. right) $\cala$-module
$\Omega^1_u(\cala)$ is isomorphic to $\cala\otimes d_u\cala$ (resp.
$d_u\cala\otimes \cala$) whereas the kernel of $d_u$ is $\Bbb C\bbbone$, i.e.
$d_u\cala \simeq \cala/\c\bbbone$ as vector space. One has the following
structure theorem for first order operators.

\begin{theorem}
Let $\calm$ and $\caln$ be two $(\cala,\calb)$-bimodules and let $D$ be a first
order operator of $\calm$ into $\caln$. Then, there is a unique
$(\cala,\calb)$-bimodule homomorphism $\sigma_L(D)$ of $\Omega^1_u(\cala)
\otimesinf_{\cala}\calm$ into $\caln$ and there is a unique
$(\cala,\calb)$-bimodule homomorphism $\sigma_R(D)$ of
$\calm\otimesinf_{\calb}\Omega^1_u(\calb)$ into $\caln$ such that one has:
\[
D(fmg)=fD(m)g + \sigma_L(D)(d_uf\otimes m)g + f\sigma_R(D)(m\otimes d_ug)
\]
for any $m\in \calm$, $f\in\cala$ and $g\in \calb$.
\end{theorem}

\noindent {\bf Proof}.  By definition, one has
\[
D(fmg)-D(fm)g-fD(mg)+fD(m)g=0
\]
 which can be rewritten
\[
D(fmg)=fD(m)g + (D(fm)-fD(m))g + f(D(mg)-D(m)g).\]
Now we know (lemma 1, a), etc.) that $m\mapsto D(fm)-fD(m)$ is a right
$\calb$-module homomorphism, furthermore it vanishes whenever $f\in \Bbb C
\bbbone$, so one defines a $(\cala,\calb)$-bimodule homomorphism
$\tilde\sigma_L(D)$ of $\Omega^1_u(\cala)\otimes \calm$ into $\caln$ by setting
$\tilde\sigma_L(D)(f_0d_uf_1\otimes m)=f_0(D(f_1m)-f_1D(m))$. Moreover, one has
\[
\begin{array}{l}
\tilde\sigma_L(D)(f_0d_u(f_1) h\otimes m)=\tilde\sigma_L(D)
(f_0d_u(f_1h)\otimes m)-\tilde \sigma_L(D)(f_0f_1d_uh\otimes m)\\
= f_0(D(f_1hm)-f_1h D(m))-f_0f_1(D(hm)-hD(m))\\
= f_0(D(f_1hm)-f_1D(hm))\\
= \tilde\sigma_L (D)(f_0d_uf_1\otimes hm).
\end{array}
\]
This means that $\tilde\sigma_L(D)(\alpha h\otimes
m)-\tilde\sigma_L(D)(\alpha\otimes h m)=0$, for any $\alpha\in
\Omega^1_u(\cala)$, $m\in\calm$ and $h\in\cala$, i.e. that $\tilde \sigma_L(D)$
passes to the quotient and defines a $(\cala,\calb)$-bimodule homomorphism
$\sigma_L(D)$ of $\Omega^1_u(\cala)\otimesinf_{\cala}\calm$ into $\caln$. The
uniqueness of $\sigma_L(D)$ is obvious by setting $g=\bbbone$ in the statement.
One proceeds similarily on the other side for $\sigma_R(D)$. $\square$\\

\noindent In view of 3.3, the converse is also true: any element $D$ of
$\call(\calm,\caln)$ for which there are $\sigma_L(D)$ and $\sigma_R(D)$ as
above is a first order operator.\\

\noindent It is clear that $\sigma_L(D)$ and $\sigma_R(D)$ are the appropriate
generalization of the notion of symbol in this setting. We shall refer to them
as the {\sl left} and the {\sl right universal symbol of} $D$ respectively.\\

\noindent In order to make contact with this notion, let us investigate the
examples of the last section. So let $D$ be a first order operator of $\cala$
into a bimodule $\caln$ and let $p_R(D)$ and $p_L(D)$ be the corresponding
derivations of $\cala$ into $\caln$ as in 3.1, then $\sigma_D(D)=i_{p_R(D)}$
and $\sigma_L(D)=i_{p_L(D)}$ are the canonical homomorphisms of bimodule over
$\cala$ of $\Omega^1_u(\cala)$ into $\caln$ associated with the derivations
$p_R(D)$ and $p_L(D)$, (i.e. such that $p_R(D)=i_{p_R(D)}\circ d_u$ and
$p_L(D)=i_{p_L(D)}\circ d_u$); in particular, if $\delta$ is a derivation of
$\cala$ into $\caln$, one has $\sigma_R(\delta)=\sigma_L(\delta)=i_\delta$. In
the cases of the examples of 3.3, for the case 1, one has
$\sigma_L(D)=\sigma_L\circ(i_{d_L}\otimes id_{\calm})$ and the case 2, one has
$\sigma_R(D)=\sigma_R\circ (id_{\calm}\otimes i_{d_R})$, where $i_{d_L}$ is the
canonical homomorphism of bimodule over $\cala$ of $\Omega^1_u(\cala)$ into
$\Omega^1_L$ associated with the derivation $d_L$, $i_{d_R}$ is the canonical
homomorphism of bimodule over $\calb$ of $\Omega^1_u(\calb)$ into $\Omega^1_R$
associated with the derivation $d_R$ and where $id_{\calm}$ is the identity
mapping of $\calm$ onto itself.\\

\noindent {\bf Remark}\\
In the case where the right and the left module structures are related, the
$\sigma_R(D)$ and the $\sigma_L(D)$ are also related. For instance if
$\cala=\calb$ is a commutative algebra and if $\calm$ and $\caln$ are
$\cala$-modules then $\sigma_R(D)=\sigma_L(D)\circ T$ where $T$ is the
transposition of $\calm\otimesinf_{\cala}\Omega^1_u(\cala)$ onto
$\Omega^1_u(\cala)\otimesinf_{\cala}\calm$ defined by  $T(m\otimes
\alpha)=\alpha \otimes m$, $\forall m\in \calm$ and $\forall
\alpha\in\Omega^1_u(\cala)$. In particular, this applies if
$\cala=\calb=C^\infty(X)$ and if $\calm$ and $\caln$ are the smooth sections of
smooth vector bundles $E$ and $F$ over the smooth manifold $X$; in this case
one has $\sigma_R(D)=\sigma(D)\circ (i_d\otimes id_{\calm})$ where $\sigma(D)$
is the usual symbol of the first order linear partial differential operator $D$
of $E$ into $F$ and $i_d$ is the canonical homomorphism of
$\Omega^1_u(C^\infty(X))$ into $\Omega^1(X)$ (= the ordinary 1-forms on $X$)
associated with the ordinary differential $d:C^\infty(X)\rightarrow
\Omega^1(X)$.

\section{Application to the theory of connections}

Let $\Omega$ be a graded differential algebra, with differential $d$, such that
$\Omega^0=\cala$. The restriction of the differential $d$ to $\cala$ is then a
derivation of $\cala$ into the bimodule $\Omega^1$ over $\cala$ (the
$\Omega^n$, $n\in\Bbb  N$, are bimodules over $\cala$). Recall that a
$\Omega$-connection on a left $\cala$-module $\calm$ [3], [4] is a linear
mapping $\nabla$ of $\calm$ into $\Omega^1\otimesinf_{\cala}\calm$ satisfying
\[
\nabla(fm)=f\nabla(m)+df\otimes m,\ \forall m\in \calm\ \mbox{and}\ \forall
f\in \cala.
\]
Suppose now that $\calm$ is not only a left $\cala$-module but is a bimodule
over $\cala$. It follows then from 3.3 that a left $\cala$-module
$\Omega$-connection $\nabla$ on $\calm$ as above is a first order operator of
the bimodule $\calm$ into the bimodule $\Omega^1\otimesinf_{\cala}\calm$ over
$\cala$. Consequently, in view of the structure theorem of 4, there is a
homomorphism of bimodule over $\cala$, $\sigma_R(\nabla)$, of
$\calm\otimesinf_{\cala}\Omega^1_u(\cala)$ into
$\Omega^1\otimesinf_{\cala}\calm$ such that
\[
\nabla(mg)=\nabla(m)g + \sigma_R(\nabla)(m\otimes d_ug),\ \forall m\in\calm\
\mbox{and}\ \forall g\in \cala.
\]
The homomorphism $\sigma_R(\nabla)$ is the right universal symbol of $\nabla$
whereas its left universal symbol is simply $\sigma_L(\nabla)=i_d\otimes
id_{\calm}$ where $i_d$ is the canonical bimodule homomorphism of
$\Omega^1_u(\cala)$ into $\Omega^1$ induced by $d$, (i.e. such that $d=i_d\circ
d_u$), and $id_{\calm}$ is the identity mapping of $\calm$ onto itself.\\
In the case where $\sigma_R(\nabla)$ factorizes through a bimodule homomorphism
$\sigma$ of $\calm\otimesinf_{\cala}\Omega^1$ into
$\Omega^1\otimesinf_{\cala}\calm$ as
$\sigma_R(\nabla)=\sigma\circ(id_{\calm}\otimes i_d)$, we call $\nabla$ {\sl a
(left) bimodule} $\Omega$-{\sl connection on the bimodule} $\calm$. (One
defines similarily a (right) bimodule $\Omega$-connection by starting with a
right $\cala$-module $\Omega$-connection on $\calm$).\\
In the case where $\calm$ is the bimodule $\Omega^1$ itself, a (left) bimodule
$\Omega$-connection on $\Omega^1$ in the above sense is just what is the first
part of the proposal of J. Mourad [9] for the definition of linear connections
in noncommutative geometry, (the second part which relates $\sigma$ and the
product $\Omega^1\otimesinf_{\cala}\Omega^1\rightarrow \Omega^2$  makes sense
only in this case). This proposal was applied in some simple examples, e.g.
[5], [8].\\
On the other hand, in [7] connections on central bimodules for the
derivation-based differential calculus i.e. for $\Omega=\Omega_{\der}(\cala)$
were defined and it was pointed out in [9] that when the central bimodule is
$\Omega^1_{\der}(\cala)$ itself the connections of [7] are linear connections
in the sense of [9], (the choice of $\sigma$ being fixed). By a similar
argument one sees that the derivation-based connections on central bimodules of
[7] are bimodules $\Omega_{\der}(\cala)$-connections in the above sense. To be
precise, this holds strictly in the finite dimensional cases, otherwise, one
has to introduce completions of the
$\Omega^1_{\der}(\cala)\otimesinf_{\cala}\calm$, but this is merely a
technicality.\\
Thus the above definition of (left) bimodule $\Omega^1$-connections seems very
general. However, it is worth noticing here that for the derivation-based
connections on central bimodules, the curvatures of these connections are
bimodule homomorphisms whereas this is generally not the case for the curvature
of a (left) bimodule $\Omega$-connection, (it is of course always a left
$\cala$-module homomorphism), although differences of such connections for a
fixed $\sigma$ are bimodule homomorphisms, (see e.g. [5]).\\

\section*{Acknowledgements}

It is a pleasure to thank John Madore for stimulating discussions.


\newpage

\begin{thebibliography}{99}


\bibitem [1]{1} N. Bourbaki.  \newblock Alg\`ebre I, Chapitre III. \newblock
Paris, Hermann 1970.

\bibitem [2]{2} H. Cartan, S. Eilenberg. \newblock Homological algebra.
\newblock Princeton University Press 1973.


\bibitem [3]{3} A. Connes.  \newblock Non commutative geometry.
\newblock Academic Press 1994.


\bibitem [4]{4} A. Connes.  \newblock Non-commutative differential
geometry.  \newblock {\em Publ. IHES} 62: 257, 1986.


\bibitem [5]{5} M. Dubois-Violette, J. Madore, T. Masson, J. Mourad.
\newblock Linear connections on the quantum plane.  \newblock{\em
Preprint L.P.T.H.E.-ORSAY} 94/94, to appear in {\em Lett. Math. Phys.}

\bibitem [6]{6} M. Dubois-Violette, P.W. Michor.  \newblock D\'erivations
et calcul diff\'erentiel non commutatif II.  \newblock {\em C.R. Acad.
Sci. Paris} 319, S\'erie I: 927-931, 1994.

\bibitem [7]{7} M. Dubois-Violette, P.W. Michor. \newblock Connections on
central bimodules. \newblock Preprint LPTHE-ORSAY 94/100 and ESI-preprint 210.

\bibitem [8]{8} J. Madore, T. Masson, J. Mourad. \newblock Linear
connections on matrix geometries.  \newblock {\em Class. and Quant. Grav.} 12,
1429 (1995).

\bibitem [9]{9} J. Mourad. \newblock Linear connections in
noncommutative geometry.  \newblock  {\em
Class. and Quant. Grav.} 12, 965-974 (1995).



\end{thebibliography}
\end{document}